\documentclass[12pt]{article}
\usepackage{amsfonts}
\usepackage{mathrsfs}
\usepackage{amsmath}
\usepackage{latexsym}
\usepackage{amssymb}
\usepackage{amsthm}
\usepackage{amscd}

\title{On Pointed Hopf Algebras with Sporadic Simple Groups ${\rm HS}$ and ${\rm Co3}$}
\author{  \small Shouchuan Zhang,    \ \ Jing Cheng
\\
\small  Department  of Mathematics, Hunan University\\  \small
Changsha  410082, \ P.R. China
 }
\date{}

\begin{document}
\newtheorem{Proposition}{\quad Proposition}[section]
\newtheorem{Theorem}{\quad Theorem}
\newtheorem{Definition}[Proposition]{\quad Definition}
\newtheorem{Corollary}[Proposition]{\quad Corollary}
\newtheorem{Lemma}[Proposition]{\quad Lemma}
\newtheorem{Example}[Proposition]{\quad Example}
\newtheorem{Remark}[Proposition]{\quad Remark}

\maketitle \addtocounter{section}{-1}

\numberwithin{equation}{section}

\date{}

\begin {abstract} Every non quasi- $-1$-type Nichols algebra
 is infinite dimensional. All quasi- $-1$-type Nichols algebra over
sporadic simple groups ${\rm HS}$ and ${\rm Co3}$ are found.


keywords: Quiver, Hopf algebra, Weyl group.
\end {abstract}

\section{Introduction}\label {s0}

This article is to contribute to the classification of
finite-dimensional complex pointed Hopf algebras  with sporadic
simple group $G= {\rm HS}$ or ${\rm Co3}$

 Many papers are about the classification of finite dimensional
pointed Hopf algebras, for example,  \cite{ AS02, AS00, AS05, He06,
AHS08, AG03, AFZ08, AZ07,Gr00,  Fa07, AF06, AF07, ZZC04, ZCZ08,
ZZWCY08}.

 Sporadic simple groups are important in  the theories of Lie groups,
Lie algebras, vertex operator algebras  and algebraic groups. Every
non quasi- $-1$-type Nichols algebra
 is infinite dimensional (see \cite [Lemma 1.8, 1.9] {AF07}). In
this paper we find all quasi- $-1$-type Nichols algebra over
sporadic simple groups ${\rm HS}$ and ${\rm Co3}$ by using the
results in \cite {ZCH09}.

For  $s\in G$ and  $(\rho, V) \in  \widehat {G^s}$, here is a
precise description of the Yetter-Drinfeld ( {\rm YD} in short )
module $M({\mathcal O}_s, \rho)$, introduced in \cite {Gr00, AZ07}.
Let $t_1 = s$, \dots, $t_{m}$ be a numeration of ${\mathcal O}_s$,
which is a conjugacy class containing $s$,  and let $g_i\in G$ such
that $g_i \rhd s := g_i s g_i^{-1} = t_i$ for all $1\le i \le m$.
Then $M({\mathcal O}_s, \rho) = \oplus_{1\le i \le m}g_i\otimes V$.
Let $g_iv := g_i\otimes v \in M({\mathcal O}_s,\rho)$, $1\le i \le
m$, $v\in V$. If $v\in V$ and $1\le i \le m $, then the action of
$h\in G$ and the coaction are given by
\begin {eqnarray} \label {e0.11}
\delta(g_iv) = t_i\otimes g_iv, \qquad h\cdot (g_iv) =
g_j(\gamma\cdot v), \end {eqnarray}
 where $hg_i = g_j\gamma$, for
some $1\le j \le m$ and $\gamma\in G^s$. The explicit formula for
the braiding is then given by
\begin{equation} \label{yd-braiding}
c(g_iv\otimes g_jw) = t_i\cdot(g_jw)\otimes g_iv =
g_{j'}(\gamma\cdot v)\otimes g_iv\end{equation} for any $1\le i,j\le
m$, $v,w\in V$, where $t_ig_j = g_{j'}\gamma$ for unique $j'$, $1\le
j' \le m$ and $\gamma \in G^s$. Let $\mathfrak{B} ({\mathcal O}_s,
\rho )$ denote $\mathfrak{B} (M ({\mathcal O}_s, \rho ))$, which is
called a bi-one type Nichols algebra (see \cite {ZZWCY08}).
$M({\mathcal O}_s, \rho )$ is a simple {\rm YD} module (see \cite
[Section 1.2 ] {AZ07}). Furthermore, if $\chi$ is the character of
$\rho$, then we also denote $\mathfrak{B} ({\mathcal O}_s, \rho )$
by $\mathfrak{B} ({\mathcal O}_s, \chi )$.

\section {Tables about  $-1$- type} \label {s6}
In this section all $-1$- type bi-one Nichols algebras over ${\rm
HS}$ and ${\rm Co3}$
 up to  graded pull-push
{\rm YD} Hopf algebra isomorphisms,  are listed in table.

\subsection {Quasi-real elements}

\begin{Definition} \label {1.1}
Let $G$ be a finite group with $s\in G$. $s$ is called a quasi-real
element if $s^2=1$ or there exists an integer number $j$ such that
$s^j \in {\mathcal O}_s$ with $s^j \not= s$. Furthermore,  $s$ is
called a strongly quasi-real element if $s^2=1$ or
 there exists an integer number $j$ such that $s^j $, $s^{j^2}$ and $  s$ are different each
other with $s^j \in {\mathcal O}_s$.

\end{Definition}

A Group is said to be quasi-real if its every element is quasi-real.
It is clear that every real element is a quasi-real element.

\begin{Definition} \label {1.2}
Let $G$ be a finite group with $s\in G$ and $s$   a quasi-real
element in $G$. $\mathfrak B({\mathcal O}_s, \rho)$ with $\rho (s) =
q_{ss} {\rm id}$ is called to be of quasi- $-1$ type if one of the
following conditions holds.

{\rm (i)} $s$ is a strongly real element and the order of $s$ is
even with $q_{ss} =-1$.

{\rm (ii)} $s$ has even order  and ${\rm deg }(\rho)>1$  with
$q_{ss}  = -1$.

{\rm (iii)}  ${\rm deg }(\rho)=1$ and the order of $s$ is even with
$q_{ss} =-1$

(iii)  ${\rm deg }(\rho)=1$ and  $q_{ss}$ is a primitive $3$-th root
of unity.
\end{Definition}

By {\rm GAP} we have the following Lemmas.

\begin{Lemma} \label {1.3}
 $s_i^j $, $s_i^{j^2}$ and $  s_i$ are
different each other with $s_i^j \in {\mathcal O}_{s_i}$ in the
following cases:

{\rm (i)} when  $G$ is  ${\rm HS}$ with   $(i,   j)$ $=$      $(8,
3)$, $(9,  3 )$,   $(10,  2 )$,   $(12,   3)$, $(13, 3 )$,  $(14,
3)$, $(15,   2)$,  $(17,   2)$,  $(22,   2)$, $(23, 3 )$,  $(24,  3
)$.

{ \rm (ii)} when $G$ is  ${\rm Co3}$ with $(i, j)=$ $(2,  2 )$, $(3,
2)$,
 $(4,  2)$,   $(5,  2 )$,   $(7,
5)$,   $(8,   2)$,
 $(14,  3 )$,
 $(15,   3)$,
 $(16,  3 )$,
 $(17,  2 )$,
 $(19,   2)$,
$(20,   2)$,   $(22,  3 )$,   $(23,  3)$,
 $(33,   3)$,
 $(34,   3)$,
 $(35,   3)$,
 $(36,   3)$,
$(38, 2 )$,
 $(39,  7)$, $(42, 2)$.

\end{Lemma}

\begin{Lemma} \label {1.4}

{\rm (i)} If   $G$ is  ${\rm HS}$,  then $s_i$ is a strongly
quasi-real element except $i =$ $2,$ $3$,  $5$, $6, $ $11,  16$,
$18$, $19,$ $20,$ $21$.

{\rm (ii)} If   $G$ is  ${\rm Co3}$,   then $s_i$ is a strongly
quasi-real element except $i =6,$ $ 9, 10,$ $  12, $ $ 18, $ $ 21,
24, $ $ 25, 26, 27, $ $  28, 29, 30,$ $  31, 32, $ $ 37, 40, 41$.

\end{Lemma}

\begin{Lemma} \label {1.5} $s_i^j $ and $  s_i$ are
different each other with $s_i^j \in {\mathcal O}_{s_i}$ in the
following cases:

(i) when  $G$ is  ${\rm HS}$ with  $(i, j)$ $=$ $(2,  3 )$, $(3,
3)$, $(5, 3)$,
 $(6,   3)$,
 $(8,  3 )$,
$(9,   2)$,   $(10,  2 )$,
 $(11,  2 )$,
 $(12,   3)$,
$(13,  3 )$,   $(14,  3 )$,   $(15,  2 )$,
 $(16,  3 )$,
 $(17,   2)$,
 $(18,   5)$,
 $(19,  5 )$,
 $(20,   5)$,
$(21,  3 )$,
 $(22,  2 )$,
 $(23,   3)$,
  $(24,   3)$.

{ \rm (ii)} when $G$ is  ${\rm Co3}$ with $(i, j)=$ $(2,  2 )$, $(3,
2)$,
 $(4,  2)$,   $(5,  2 )$,  $(6,  2 )$,  $(7,
5)$,   $(8,   2)$, $(9,  2 )$, $(10,  5 )$, $(12,  2 )$,
 $(14,  3 )$,
 $(15,   3)$,
 $(16,  3 )$,
 $(17,  2 )$, $(18,  3 )$,
 $(19,   2)$,
$(20,   2)$,   $(21,  2 )$, $(22,  3 )$,   $(23,  3)$, $(24,  3 )$,
$(25,  3 )$, $(26,  3 )$, $(27,  5 )$, $(28,  5 )$, $(29,  5 )$,
$(30, 5 )$, $(31,  3 )$, $(32,  5 )$,
 $(33,   3)$,
 $(34,   3)$,
 $(35,   3)$,
 $(36,   3)$, $(37,  5 )$,
$(38, 2 )$,
 $(39,  7)$,  $(40,  5 )$,    $(41,  5 )$, $(42, 2)$.

\end{Lemma}

\begin{Lemma} \label {1.6}  ${\rm HS}$ and  ${\rm Co3}$ are quasi-real.

\end{Lemma}

\subsection {Quasi-$-1$ types}

In  Table 1--3,   we use the following notations. $s_i$ denotes the
representative of $i$-th conjugacy class of $G$ ($G$ is  ${\rm HS}$
and ${\rm Co3}$); $\chi _i ^{(j)}$  denotes the $j$-th character of
 $G^{s_i}$ for any $i$;
  $\nu_i ^{(1)}$ denotes    the number of
conjugacy classes of  the centralizer  $G^{s_i}$; $\nu_i ^{(2)}$
denote the number of character $\chi _i^{(j)}$
 of
$G^{s_i}$ with non  $-1$-type $\mathfrak{B}(\mathcal{O}_{s_i},  \chi
_i ^{(j)})$;  ${\rm cl}_i [j]$ denote that  $s_i$ is in  $j$-th
conjugacy class of $G^{s_i}$.

Using the results in \cite {ZCH09} we can obtain the following
tables.

\begin{tabular}{|l|l|l|l|l|l|}
  \hline
  ${\rm HS}$ &  &  & & &\\\hline
  $s_{i}$ & ${\rm cl}_{i}[p]$ & $\hbox {Order}(s_{i})$ & $j$ such that $\mathfrak
B({\mathcal O}_{s_i},\chi_i^{(j)})$ is of quasi- $-1$ type
   & $\nu _i^{(1)}$ & $\nu_i^{(2)}$   \\\hline
   $s_{1}$ & ${\rm cl}_{1}[1]$ & 1 & & 24 &  24 \\\hline
  $s_{2}$ & ${\rm cl}_{2}[8]$ & 8 & 3,4 & 16 & 14 \\\hline
  $s_{3}$ & ${\rm cl}_{3}[21]$ & 4 & 9,10,11,12,13,14,15,16,17,18 & 22 & 12 \\\hline
  $s_{4}$ & ${\rm cl}_{4}[2]$ & 2 &12, 24, 25, 27, 28 & 28 & 23  \\\hline
  $s_{5}$ & ${\rm cl}_{5}[11]$ & 8 & 2,4 & 16 & 14 \\\hline
  $s_{6}$ & ${\rm cl}_{6}[30]$ & 4 & 5,6,7,8,13,14,15,16,17,20,24,25,26 & 34 & 21 \\\hline
  $s_{7}$ & ${\rm cl}_{7}[2]$ & 2 &2, 3, 4, 5. 9, 10, 11, 12, 19, 20, 23, 24, 26  & 26 & 13 \\\hline
  $s_{8}$ & ${\rm cl}_{8}[3]$ & 10 & 2,4 & 20 & 18 \\\hline
  $s_{9}$ & ${\rm cl}_{9}[23]$ & 5 &  & 25 & 25  \\\hline
  $s_{10}$ & ${\rm cl}_{10}[7]$ & 15 & & 15 & 15  \\\hline
  $s_{11}$ & ${\rm cl}_{11}[21]$ & 3 & 3,4,5,6 & 21 & 17  \\\hline
  $s_{12}$ & ${\rm cl}_{12}[4]$ & 20 & 2 & 20 & 19  \\\hline
  $s_{13}$ & ${\rm cl}_{13}[2]$ & 20 & 2 & 20 & 19 \\\hline
  $s_{14}$ & ${\rm cl}_{14}[5]$ & 10 & 3,4 & 20 & 18 \\\hline
  $s_{15}$ & ${\rm cl}_{15}[25]$ & 5 & & 26 & 26 \\\hline
  $s_{16}$ & ${\rm cl}_{16}[4]$ & 4 & 8,9,10,11,16,17,24 & 26 & 19 \\\hline
  $s_{17}$ & ${\rm cl}_{17}[4]$ & 7 & & 7 & 7 \\\hline
  $s_{18}$ & ${\rm cl}_{18}[10]$ & 12 & 2,3,5 & 12 & 9 \\\hline
  $s_{19}$ & ${\rm cl}_{19}[10]$ & 6 & 5,6,7,8,9,10,11,12,13 & 15 & 6 \\\hline
  $s_{20}$ & ${\rm cl}_{20}[16]$ & 6 & 2,4,7,8,11,12,13 & 18 & 11 \\\hline
  $s_{21}$ & ${\rm cl}_{21}[2]$ & 8 & 2,3 & 16 & 14 \\\hline
  $s_{22}$ & ${\rm cl}_{22}[17]$ & 5 & & 25 & 25 \\\hline
  $s_{23}$ & ${\rm cl}_{23}[4]$ & 11 & & 11 & 11  \\\hline
  $s_{24}$ & ${\rm cl}_{24}[8]$ & 11 & & 11 & 11 \\\hline
\end{tabular}
$$\hbox {Table } 1$$

\begin{tabular}{|l|l|l|l|l|l|}
  \hline
   ${\rm CO_3}$ &  & & & &\\\hline
  $s_{i}$ & ${\rm cl}_{i}[p]$ & $\hbox {Order}(s_{i})$ &$j$ such that $\mathfrak
B({\mathcal O}_{s_i},\chi_i^{(j)})$ is of quasi- $-1$ type
   & $\nu _i^{(1)}$ & $\nu_i^{(2)}$ \\\hline
  $s_{1}$ & ${\rm cl}_{1}[1]$ & 1 &  & 42 & 42 \\\hline
  $s_{2}$ & ${\rm cl}_{2}[4]$ & 23 &  & 23 & 23 \\\hline
  $s_{3}$ & ${\rm cl}_{3}[14]$ & 23 &  & 23 & 23 \\\hline
  $s_{4}$ & ${\rm cl}_{4}[2]$ & 15 &  & 15 & 15 \\\hline
  $s_{5}$ & ${\rm cl}_{5}[25]$ & 5 &  & 25 & 25 \\\hline
  $s_{6}$ & ${\rm cl}_{6}[3]$ & 3 &   & 48  & 48  \\\hline
  $s_{7}$ & ${\rm cl}_{7}[7]$ & 18 & 2 & 18 & 17 \\\hline
  $s_{8}$ & ${\rm cl}_{8}[23]$ & 9 &  & 30 & 30 \\\hline
  $s_{9}$ & ${\rm cl}_{9}[3]$ & 3 &   & 60 & 60 \\\hline
  $s_{10}$ & ${\rm cl}_{10}[4]$ & 6 & 2,3,4,8,11,12,14,40,43,44  & 48 & 38 \\\hline
$s_{11}$ & ${\rm cl}_{11}[2]$ & 2 & 3, 10, 12, 13, 19, 20, 22, 24,
33, 39, 41, 42, 43 & 43 & 30
\\\hline
  $s_{12}$ & ${\rm cl}_{12}[44]$ & 6 &  5,6,7,8,9,10,11,12,15,26,27,36,37,41 & 45 & 31 \\\hline
  $s_{13}$ & ${\rm cl}_{13}[2]$ & 2 &2, 3, 4, 7, 8, 11, 13, 16,
  17, 18, 22, 24, 26, 28, 30  & 30 & 15 \\\hline
  $s_{14}$ & ${\rm cl}_{14}[10]$ & 20 &  2 & 20 & 19 \\\hline
  $s_{15}$ & ${\rm cl}_{15}[14]$ & 20 & 2 & 20 & 19 \\\hline
  $s_{16}$ & ${\rm cl}_{16}[12]$ & 10 &  3,4,21 & 30 & 27 \\\hline
  $s_{17}$ & ${\rm cl}_{17}[2]$ & 5 &   & 32 & 32 \\\hline
  $s_{18}$ & ${\rm cl}_{18}[4]$ & 4 &  6,11,12,13,14,17,18,23,25,32,33 & 37 & 26 \\\hline
  $s_{19}$ & ${\rm cl}_{19}[5]$ & 21 &  & 21 & 21 \\\hline
  $s_{20}$ & ${\rm cl}_{20}[19]$ & 7 &  & 21 & 21 \\\hline
 $s_{21}$ & ${\rm cl}_{21}[3]$ & 3 & 4,5,6,7,8,9 & 33 & 27 \\\hline
  $s_{22}$ & ${\rm cl}_{22}[4]$ & 14 & 2 & 14 & 13 \\\hline
  $s_{23}$ & ${\rm cl}_{23}[19]$ & 10 & 2,4 & 20 & 18 \\\hline
  $s_{24}$ & ${\rm cl}_{24}[2]$ & 8 &  2,3,4,5 & 20 & 16 \\\hline
$s_{25}$ & ${\rm cl}_{25}[3]$ & 4 & 5,6,7,8,11,12,17,18,19,20,25,26,
& 50 & 30 \\&&&27,28,35,36,39,40,43,44 & & \\\hline
 $s_{26}$ & ${\rm cl}_{26}[29]$ & 8 & 3,4,9,21,24 & 32 & 27 \\\hline
\end{tabular}
$$\hbox {Table } 2$$

\begin{tabular}{|l|l|l|l|l|l|}
  \hline
   ${\rm CO_3}$ &  & & & &\\\hline
  $s_{i}$ & ${\rm cl}_{i}[p]$ & $\hbox {Order}(s_{i})$ & $j$ such that $\mathfrak
B(\mathcal{O}_{s_i},\chi_i^{(j)})$ is of quasi- $-1$ type
   & $\nu _i^{(1)}$ & $\nu_i^{(2)}$ \\\hline

  $s_{27}$ & ${\rm cl}_{27}[16]$ & 24 &  2,3,5 & 24 & 21 \\\hline
  $s_{28}$ & ${\rm cl}_{28}[21]$ & 12 & 5,6,7,8,17,18,19,20,21,22,23,24 & 30 & 18\\\hline
  $s_{29}$ & ${\rm cl}_{29}[6]$ & 6 &  3,4,5,6,7,8,9,10,44,49 & 51 & 41 \\\hline
  $s_{30}$ & ${\rm cl}_{30}[20]$ & 24 & 2,3,5 & 24 & 21 \\\hline
  $s_{31}$ & ${\rm cl}_{31}[32]$ & 8 & 2,4,10,22,23 & 32 & 27 \\\hline
  $s_{32}$ & ${\rm cl}_{32}[17]$ & 6 & 2,4,6,8,11,12 & 20 & 14 \\\hline
  $s_{33}$ & ${\rm cl}_{33}[13]$ & 22 & 2 & 22 & 21 \\\hline
  $s_{34}$ & ${\rm cl}_{34}[8]$ & 22 & 2 & 22 & 21 \\\hline
  $s_{35}$ & ${\rm cl}_{35}[21]$ & 11 &   & 22 & 22 \\\hline
  $s_{36}$ & ${\rm cl}_{36}[12]$ & 11 &  & 22 & 22 \\\hline
  $s_{37}$ & ${\rm cl}_{37}[17]$ & 12 & 2,5,6,9,10,13,14,15,16,37 & 42 & 32 \\\hline
  $s_{38}$ & ${\rm cl}_{38}[19]$ & 15 &  & 30 & 30 \\\hline
  $s_{39}$ & ${\rm cl}_{39}[29]$ & 30 &  2 & 30 & 29 \\\hline
  $s_{40}$ & ${\rm cl}_{40}[24]$ & 12 &  2,5,6,7,8,15,16,17,18 & 36 & 27 \\\hline
  $s_{41}$ & ${\rm cl}_{41}[23]$ & 6 & 2,5,6,9,10,13,14,15,16,19 & 24 & 14 \\\hline
  $s_{42}$ & ${\rm cl}_{42}[33]$ & 9 &  & 33 & 33 \\\hline
\end{tabular}
$$\hbox {Table } 3$$

Table 1 is called the table of ${\rm HS}$ and Table 2--3 are called
the tables of ${\rm Co3}$.

\section{ Bi-one Nichols algebras over ${\rm HS}$ and ${\rm Co3}$
}\label {s6} In this section all  $-1$-type bi-one Nichols algebra
over
 ${\rm HS}$ and ${\rm Co3}$ of exceptional type up to  graded pull-push
{\rm YD} Hopf algebra isomorphisms are  given.

\begin {Lemma}\label {2.1} Assume that $s\in G$ is quasi-real with
$\rho \in \widehat {G^s}$. If $ \mathfrak{B}(\mathcal{O}_{s},\rho)$
is not  of quasi-$-1$ type, then ${\rm dim}
\mathfrak{B}(\mathcal{O}_{s},\rho)=\infty$.
\end {Lemma}

{\bf Proof}.  If follows from \cite [Lemma 1.8, 1.9] {AF07}. $\Box$.

We give our main result.

\begin {Theorem} \label {1} Let $G$ be one of ${\rm HS}$ and  ${\rm Co3}$.

 {\rm (i)}
$\mathfrak{B} ({\mathcal O}_{s_i},   \chi _i^{(j)})$ is of $-1$-type
if and only if $j$ appears in the  fourth column of the table of
$G$.

{\rm (ii)} ${\rm dim } (\mathfrak{B} ({\mathcal O}_{s_i},   \chi
_i^{(j)})) = \infty $ if $j$  does not appears in the fourth column
of the table of $G$.

\end {Theorem}

{\bf Proof.} {\rm (i)} It follows from the program.

{\rm (ii)} It follows from Lemma \ref {2.1}.
 $\Box$

\section {Appendix}

In this section Suzuki group ${\rm Sz(8)}$ is considered.

\begin{Lemma} \label {3.1}

{\rm (i)} $s^j $, $s^{j^2}$ and $  s$ are different each other with
$s^j \in {\mathcal O}_s$ in the following cases:

 $(i,   j)$ $=$   $(5, 5)$,
 $(6,  5)$,
  $(7,  5)$,
  $(11,  2)$.

{\rm (ii)} $s_i$ is strongly quasi-real except $i= 2, 3, 4, 8, 9.$

  (iii) $s_i$ is strongly quasi-real if and only if $s_i$ is
  quasi-real.

\end{Lemma}

\begin{tabular}{|l|l|l|l|l|l|}
  \hline
  ${\rm Sz(8)}$ &  &  & & &\\\hline
  $s_{i}$ & ${\rm cl}_{i}[p]$ & $Order(s_{i})$ & $j$ such that $\mathfrak
B(\mathcal{O}_{s_i},\chi_i^{(j)})$ is of quasi- $-1$ type
   & $\nu _i^{(1)}$ & $\nu_i^{(2)}$   \\\hline
   $s_{1}$ & ${\rm cl}_{1}[1]$ & 1 & & 11 & 11 \\\hline
  $s_{2}$ & ${\rm cl}_{2}[8]$ & 7 &  & 7 & 7 \\\hline
  $s_{3}$ & ${\rm cl}_{3}[21]$ & 7 &  & 7 & 7 \\\hline
  $s_{4}$ & ${\rm cl}_{4}[2]$ & 7 & & 7 & 7  \\\hline
  $s_{5}$ & ${\rm cl}_{5}[11]$ & 13 & & 7 &7 \\\hline
  $s_{6}$ & ${\rm cl}_{6}[30]$ & 13 & &13 &13 \\\hline
  $s_{7}$ & ${\rm cl}_{7}[2]$ & 13 &  & 13 & 13 \\\hline
  $s_{8}$ & ${\rm cl}_{8}[3]$ & 4 &  & 16 & 16 \\\hline
  $s_{9}$ & ${\rm cl}_{9}[23]$ & 4 &  & 16 & 16  \\\hline
  $s_{10}$ & ${\rm cl}_{10}[7]$ & 2 &9, 10, 13, 14, 17, 18, 19, 20 & 22 &  14 \\\hline
  $s_{11}$ & ${\rm cl}_{11}[21]$ & 5 & & 5 & 5  \\\hline
  \end{tabular}
$$\hbox {Table } 4$$

\begin {Proposition} \label {3.2} Let $G$ be  ${\rm Sz(8)}$ with $i = 1, 5, 6, 7, 10, 11$.

 {\rm (i)}
$\mathfrak{B} ({\mathcal O}_{s_i},   \chi _i^{(j)})$ is of
quasi-$-1$-type if and only if $j$ appears in the  fourth column of
the table of Table 4.

{\rm (ii)} ${\rm dim } (\mathfrak{B} ({\mathcal O}_{s_i},   \chi
_i^{(j)})) = \infty $ if $j$  does not appears in the fourth column
of the table of Table 4.

\end {Proposition}


\vskip 1.0cm


\begin {thebibliography} {200}
\bibitem [AF06]{AF06} N. Andruskiewitsch and F. Fantino,   On pointed Hopf algebras
associated to unmixed conjugacy classes in Sn,   J. Math. Phys. {\bf
48}(2007),    033502-1-- 033502-26. Also math.QA/0608701.


\bibitem[Atlas]{Atlas} Atlas of finite group representation- Version 3,
http://brauer.maths.qmul.ac.uk/Atlas/v3.


\bibitem [AF07]{AF07} N. Andruskiewitsch,
F. Fantino,      On pointed Hopf algebras associated with
alternating and dihedral groups,   preprint,   arXiv:math/0702559.

\bibitem [AFZ]{AFZ08} N. Andruskiewitsch,
F. Fantino,     Shouchuan Zhang,    On pointed Hopf algebras
associated with symmetric  groups,   Manuscripta Mathematica,
accepted. Also arXiv:0807.2406.



\bibitem[AG03]{AG03} N. Andruskiewitsch and M. Gra\~na,
From racks to pointed Hopf algebras,   Adv. Math. {\bf 178}(2003),
177-243.


\bibitem [AHS08]{AHS08} N. Andruskiewitsch,   I. Heckenberger,   H.-J. Schneider,
  The Nichols algebra of a semisimple Yetter-Drinfeld module,    preprint,
arXiv:0803.2430.


\bibitem [AS98]{AS98} N. Andruskiewitsch and H. J. Schneider,
Lifting of quantum linear spaces and pointed Hopf algebras of order
$p^3$,    J. Alg. {\bf 209} (1998),   645--691.

\bibitem [AS02]{AS02} N. Andruskiewitsch and H. J. Schneider,   Pointed Hopf algebras,
new directions in Hopf algebras,   edited by S. Montgomery and H.J.
Schneider,   Cambradge University Press,   2002.

\bibitem [AS00]{AS00} N. Andruskiewitsch and H. J. Schneider,
Finite quantum groups and Cartan matrices,   Adv. Math. {\bf 154}
(2000),   1--45.

\bibitem[AS05]{AS05} N. Andruskiewitsch and H. J. Schneider,
On the classification of finite-dimensional pointed Hopf algebras,
 Ann. Math.,   accepted. Also   {math.QA/0502157}.

\bibitem [AZ07]{AZ07} N. Andruskiewitsch and Shouchuan Zhang,   On pointed Hopf
algebras associated to some conjugacy classes in $S_n$,   Proc.
Amer. Math. Soc. {\bf 135} (2007),   2723-2731.




\bibitem  [CR02]{CR02} C. Cibils and M. Rosso,    Hopf quivers,   J. Alg. {\bf  254}
(2002),   241-251.

\bibitem [CR97] {CR97} C. Cibils and M. Rosso,   Algebres des chemins quantiques,
Adv. Math. {\bf 125} (1997),   171--199.

\bibitem[DPR]{DPR} R. Dijkgraaf,   V. Pasquier and P. Roche,
Quasi Hopf algebras,   group cohomology and orbifold models, Nuclear
Phys. B Proc. Suppl. {\bf 18B} (1991),   pp. 60--72.


\bibitem [Fa07] {Fa07}  F. Fantino ,   On pointed Hopf algebras associated with the
Mathieu simple groups,   preprint,    arXiv:0711.3142.

\bibitem[GAP]{GAP} The {\rm GAP }-Groups, Algorithms, and Programming, Version
 4.4.12;  2008,  http://www.gap-system.org.

\bibitem[Gr00]{Gr00} M. Gra\~na,   On Nichols algebras of low dimension,
 Contemp. Math.,    {\bf 267}  (2000),  111--134.



\bibitem[He06]{He06} I. Heckenberger,   { Classification of arithmetic
root systems},   preprint,   {math.QA/0605795}.


\bibitem[HS]{HS} I. Heckenberger and H.-J. Schneider,   { Root systems
and Weyl groupoids for  Nichols algebras},   preprint
{arXiv:0807.0691}.


\bibitem [Ra]{Ra85} D. E. Radford,   The structure of Hopf algebras
with a projection,   J. Alg. {\bf 92} (1985),   322--347.

 \bibitem [Sw] {Sw69} M. E. Sweedler,   Hopf algebras,   Benjamin,   New York,   1969.



\bibitem [ZCZ]{ZCZ08} Shouchuan Zhang,    H. X. Chen and Y.-Z. Zhang,
Classification of  quiver Hopf algebras and pointed Hopf algebras of
type one,   preprint arXiv:0802.3488.

\bibitem [ZZWCY08]{ZZWCY08} Shouchuan Zhang,  Y.-Z. Zhang,  Peng Wang,   Jing Cheng,   Hui Yang,
On Pointed Hopf Algebras with Weyl Groups of Exceptional, Preprint
arXiv:0804.2602.




\bibitem [ZWCY08a]{ZWCY08a}, Shouchuan Zhang,    Peng Wang,   Jing Cheng,   Hui Yang, The character tables of centralizers in Weyl Groups of $E_6$, $E_7$,
$F_4$,   $G_2$,   Preprint  arXiv:0804.1983.

\bibitem [ZWCYb]{ZWCY08b} Shouchuan Zhang,    Peng Wang,   Jing Cheng,   Hui Yang,
The character tables of centralizers in Weyl Group of $E_8$: I - V,
Preprint. arXiv:0804.1995,   arXiv:0804.2001,   arXiv:0804.2002,
arXiv:0804.2004,   arXiv:0804.2005.

\bibitem [ZZC]{ZZC04} Shouchuan Zhang,   Y.-Z. Zhang and H. X. Chen,   Classification of PM quiver
Hopf algebras,   J. Alg. Appl. {\bf 6} (2007)(6),   919-950. Also
math.QA/0410150.


\bibitem [ZCH]{ZCH09} Shouchuan Zhang,  Jing Cheng,  Jieqiong He,
The character tables of  centralizers in  Sporadic Groups of ${\rm
HS}$ and  ${\rm CO_3}$, Preprint arXiv:0906.0408.



\end {thebibliography}
\end {document}